\documentstyle{amsart}

\newtheorem{cor}{Corollary}
\newtheorem{th}{Theorem}
\newtheorem{thh}{Theorem H}
\newtheorem{prp}{Proposition}
\newtheorem{conj}{Conjecture}

\newtheorem{lem}{Lemma}
\newcommand{\MM}{\mbox{ \bf M}}
\newcommand{\Sym}{\mbox{ \rm Sym}}
\newcommand{\chm}{\cite{chm} }
\newcommand{\rar}{\rightarrow}
\newcommand{\BP}{{\Bbb{P}^1}}

\begin{document}
\begin{center}{\bf Lang's conjectures, Conjecture H, and uniformity} \\  by Dan
Abramovich\footnote{Partially supported by NSF grant DMS-9503276.}\\ May 22,
1995$^*$
\end{center}

{\sc $^*$ The purpose of this note is to wish a happy birthday to Professor
Lucia Caporaso.}

We prove that Conjecture H of Caporaso et al. \cite{chm}, \S 6 together with
Lang's
conjecture implies the uniformity of rational points on varieties of
general type, as predicted in \cite{chm}; a few applications in arithmetic
and geometry are stated.  As in \cite{abr},  one uses the fact
that rational points on
$X^n_B$ are $n$-tuples of points when trying to bound $n$, as well as the fact
that the Noetherian induction used in \cite{chm} can be wired into the
definitions.

Geometers with number-theory anxiety should skip directly to the corollary on
the last page.

\section{Introduction}
Let $X$ be a variety of general type defined over a number field $K$. It was
conjectured by S. Lang that the set of rational points $X(K)$ is not Zariski
dense in $X$.
In the paper \cite{chm} of L. Caporaso, J. Harris and B. Mazur
 it is shown  that the above conjecture of Lang implies
the existence of a uniform bound on the number of $K$-rational points of all
curves of fixed genus $g$ over $K$.

The paper \chm has  immediately created a chasm
among arithmetic geometers.
This chasm, which often runs right in the
middle of the personalities involved, divides between loyal believers of Lang's
conjecture, who marvel in this powerful implication, and the disbelievers, who
try to use this implication to derive counterexamples to the conjecture.

In this paper we will attempt to deepen this chasm, using the
techniques of \cite{chm} and continuing \cite{abr}, by proving more
implications, some of which very strong, of various conjectures of Lang. Along
the way we will often use a conjecture donned by Caporaso et al. {\em
Conjecture H} (see again \cite{chm}, \S 6) about {\em H}igher dimensional
varieties,
which is regarded very plausible among experts of higher dimensional algebraic
geometry.

Before we state any results, we need to specify  various conjectures which we
will apply.

\subsection{A few conjectures of Lang} Let $X$ be a variety of general type
over a field $K$ of characteristic 0. In view of Faltings's proof of Mordell's
conjecture, Lang has
stated the following conjectures:

\begin{conj} \begin{enumerate} \item (Weak Lang conjecture) If $K$ is finitely
generated over $\Bbb{Q}$ then the set of rational points $X(K)$ is not Zariski
dense in $X$.
\item (Weak Lang conjecture for function fields) If $k\subset K$
is a finitely generated regular extension in characteristic 0, and if $X(K)$ is
Zariski dense in $X$, then $X$ is birational to a variety $X_0$ defined over
$k$ and the {\em ``non-constant points''} $X(K)\setminus X_0(k)$ are not
Zariski dense in $X$.
\item (Geometric Lang's conjecture) Assuming only
$Char(K) = 0$, there is a proper Zariski closed subset $Z(X) \subset X$, called
in \cite{chm} the {\em Langian exceptional set}, which is the union of all
positive dimensional subvarieties which are not of general type.
\item (Strong Lang conjecture) If $K$ is finitely generated over $\Bbb{Q}$
then there is a Zariski  closed subset $Z\subset X$ such that for any finitely
generated field $L\supset K$ we have that $X(L)\setminus Z(L)$ is finite.
\end{enumerate}
\end{conj}

These conjectures and the relationship between them are studied in
\cite{langbul}, \cite{lang3} and in the introduction of \cite{chm}. For
instance, it should be
noted that the weak Lang conjecture together with the geometric conjecture
imply the strong Lang conjecture.

It should also be remarked that the analogous situation over fields of positive
characteristic is subtle and interesting. See a recent survey by Voloch
\cite{voloch}.

\subsection{Conjecture H} An important tool used by Caporaso et al. in
\cite{chm} is that of fibered powers. Let $X\rar B$ be a morphism of varieties
in characteristic 0, where the general fiber is a variety of general type. We
denote by $X^n_B$ the
$n$-th fibered power of $X$ over $B$.

\begin{conj} (Conjecture H of \cite{chm}) For sufficiently large $n$, there
exists a dominant rational map $h_n:X^n_B \dashrightarrow W_n$ where $W_n$ is a
variety of general type, and where the restriction of $h_n$ to the general
fiber $(X_b)^n$ is generically finite.
\end{conj}

This conjecture is known for curves and surfaces:
\begin{thh} (Correlation theorem of \cite{chm}) Conjecture H holds when  $X\rar
B$ is a family of curves  of genus $>1$. \end{thh}
\begin{thh} (Correlation theorem of \cite{hassett}) Conjecture H holds when
$X\rar B$ is a family of surfaces of general type . \end{thh}

Using their  Theorem H 1, and  Lemma 1.1 of \cite{chm}, Caporaso
et al. have shown that the weak Lang conjecture implies a uniform bound on the
number of rational points on curves (Uniform bound theorem, \cite{chm} Theorem
1.1).

It should be noted that the proofs of theorems H 1 and H 2 give a bit more:
they describe a natural dominant rational map $X^n_B\rar W$.
For
the case of curves, if $B_0$ is the image of $B$ in the moduli space, $
\MM_g$, then for sufficiently large $n$ the inverse image $B_n\subset
{\MM}_{g,n}$ in the moduli space of $n$-pointed curves is a variety of general
type. Therefore the moduli map $X^n_B \dashrightarrow B_n\subset
\MM_{g,n}$ satisfies the requirements. A similar construction works for
surfaces of general type, and one may ask whether this should hold in general.

It is convenient to make the following definitions when discussing Lang's
conjectures: \\
{\bf Definition:} {\em 1. A variety $X/K$ is said to be a {\bf Lang variety} if
there is a
dominant rational map $X_{\overline{K}} \dashrightarrow W$, where $W$ is a
positive
dimensional variety of general type.

2. A positive dimensional variety $X$ is said to be {\bf geometrically
mordellic}  (In short  GeM) if $X_{\overline{K}}$ does not
contain subvarieties which are not of general type.}

In \cite{lang3}, in the course of stating even more far reaching conjectures,
Lang defined by a notion of {\em algebraically hyperbolic}
varieties which is very similar, and conjecturally the same as that of GeM
varieties. I chose to use a different terminology here, to avoid confusion.

Note that the weak Lang conjecture directly implies that the rational points on
a Lang variety over a number field are not Zariski dense, and that there are
only finitely many rational points over a number field on a GeM variety.

\subsection{Summary of results}

An indicated in \cite{chm} \S 6, Conjecture H together with Lang's conjectures
should have very strong implications for counting rational points on varieties
of general type, similar to the uniform bound theorem of \chm. Here we will
prove the following basic result:

\begin{th}\label{unif} Assume that the weak Lang conjecture as well as
conjecture H hold.
Let $X \rar B$ be a family of GeM varieties over a number field $K$ (or any
finitely generated field over $Q$). Then there is a uniform bound on $\sharp
X_b(K)$.
\end{th}

One may refine this theorem for arbitrary families of varieties of
general type, obtaining a bound on the number of points which do not lie in
lang exceptional sets of fibers. If one assumes Lang's geometric conjecture,
one obtains a closed subset $Z(X_b)$ for every $b\in B$. A natural question
which arises in such a refinement is: how do these subsets fit together? An
answer was given in \chm, Theorem 6.1, assuming conjecture H as well: the
varieties $Z(X)$ are uniformly bounded. We will
show that, using results of Viehweg, one does not need to assume conjecture H:

\begin{th}\label{Z(X)}(Compare \chm, Theorem 6.1) Assume that the geometric
Lang conjecture holds. Let $X\rar B$ be a family of varieties of general type.
Then there is a proper closed subvariety $\tilde{Z}\subset X$ such that for any
$b\in B$ we have $Z(X_b)\subset \tilde{Z}$.
\end{th}

Using theorem \ref{Z(X)}, we can apply theorem \ref{unif} to any family  $X\rar
B$ of varieties of general type, assuming that the geometric Lang conjecture
holds: we can bound the rational points in the complement of $\tilde{Z}$.

We will apply our theorem \ref{unif} in various natural cases. An immediate but
rather surprising application is the following theorem:

\begin{th}\label{unideg} Assume that the weak Lang conjecture as well as
conjecture H hold.
Let $X \rar B$ be a family of GeM varieties over a field $K$ finitely generated
over $Q$. Fix a number $d$.  Then there is a uniform bound $N_d$ such that for
any field extension $L$ of $K$ of degree $d$ and every $b\in B(L)$ we have
$\sharp X_b(L)<N_d$.
\end{th}

As a corollary, we see that Lang's conjecture together with conjecture H
imply the existence of a bound on the number of points on curves of fixed
genus $g$ over a number field $K$ which depends only on the degree of the
number field.

These results have natural analogues for function fields. We will state a few
of these, notably:

\begin{th}\label{unigon} Assume that Lang's conjecture for function fields
holds.  Fix an integer  $g>1$.
Then there is an integer $N(g)$ such that for any generically smooth
fibration of
curves $C\rar D$ where the fiber has  genus $g$ and the base is hyperelliptic
curve, there are at most $N$ non-constant sections $s:D\rar C$.
\end{th}

We remind the reader that the {\em gonality } of a curve $D$ is the minimal
degree of a nonconstant rational function on $D$ (so a curve of gonality 2 is
 hyperelliptic). One expects the above theorem to be generalized to the
situation where ``hyperelliptic curve'' is replaced by ``curve of gonality
$\leq d$'' for fixed $d$.

\section{Proof of theorem \ref{unif}}
\subsection{Preliminaries}
Throughout this subsection {\bf we assume that conjecture H holds}, and the
base field is algebraically closed.

Observe that a positive dimensional subvariety of an GeM variety is GeM; and
the normalization of an GeM variety is GeM. Note also that a variety dominating
a Lang variety is a Lang variety as well.

\begin{prp} Let $X \rar B$ be a family of GeM varieties. Let $F \subset X$ be a
reduced subscheme such that every component of $F$ dominating $B$ has positive
fiber dimension. Then for $n$
sufficiently large, every component of the fibered power $F^n_B$ which
dominates $B$ is a Lang variety.
\end{prp}

The proof will use the following lemmas:

\begin{lem} Let $X\rar B$ and $F$ be as above, and assume that the general
fiber of $F\rar B$ is irreducible. Then for $n$ sufficiently large, the
dominant component of $F^n_B$ is a Lang variety.
\end{lem}

{\bf Proof.} Apply conjecture $H$ to $F \rar B$, using the fact that the fibers
of $F$ are of general type.

\begin{lem} Let $X\rar B$ and $F$ be as in the proposition, with $F$
irreducible. Then for $n$ sufficiently large, every component of the fibered
power $F^n_B$ which dominates $B$ is a Lang variety.
\end{lem}

{\bf Proof.} Let $\tilde{F}$ be the normalization of $F$, and let $\tilde{F}
\rar \tilde{B} \rar B$ be the Stein factorization. Denote by $c$ the degree of
$\tilde{B}$ over $B$. Let $G \subset \tilde{F}^n_B$ be a dominant component.
Then $G$ parametrizes $n$-tuples of points in the fibers of $\tilde{F}$ over
$B$, and since $G$ is irreducible, there is a decomposition $\displaystyle
\{1,\ldots,n\} = \cup_{i=1}^c J_i$ and $G$ surjects onto the dominant component
of $\tilde{F}^{J_i}_{\tilde{B}}$. At least one of $J_i$ has at least $n/c$
elements. Using lemma 1 applied to $\tilde{F}\rar \tilde{B}$, we see
that for $n/c$ large enough $G$ is a Lang variety.

{\bf Proof of proposition. } Let $F= F_1 \cup\ldots \cup F_m$ be the
decomposition into irreducible components. Let $G$ be a dominant component of
$F^n_B$. Then $G$ dominates $(F_1)^{n_1}_B \times_B
\cdots\times_B(F_m)^{n_m}_B$. For at least one $i$ we have $n_i > n/m$, so
applying the previous
lemma we obtain that $G$ is a Lang variety.

\subsection{Prolongable points}
We return to the setup in  theorem \ref{unif}.

{\bf Definition.} 1. A point $x_n = (P_1,\ldots,P_n)\in X^n_B(K) $ is said to
be off diagonal if for any $1\leq i< j\leq n$ we have $P_i\neq P_j$. We extend
this for $n=0$ trivially by agreeing that any point of $B(K)$ is off diagonal.

2. Let $m>n$. An off diagonal point $x_n$ is said to be $m$-prolongable if
there is an off-diagonal $x_m\in X^m_B(K)$ whose first $n$ coordinates give
$x_n$.

Let $E_n^{(m)}$ be the set of $m$-prolongable points on $X^n_B$, and let
$F_n^{(m)}$ be the Zariski closure. Let $ F_n = \displaystyle\cap_{m>n}
F_n^{(m)}$. By the Noetherian property  of the Zariski topology we have $F_n =
F_n^{(m)}$ for some $m$.

All we need to show is $F_n = \emptyset$ for some $n$.

\begin{lem} We have a surjection $F_{n+1} \rar F_n$.
\end{lem}

{\bf Proof.} The set $E_{n+1}^{(m)}$ surjects to $E_n^{(m)}$ for any $m>n+1$.

\begin{lem} Every  fiber of $F_{n+1} \rar F_n$ is positive dimensional.
\end{lem}

{\bf Proof.} Suppose that over an open set in $F_n$ the degree of the map is
$d$. Then $E_n^{(n+d+1)}$ cannot be dense in $F_n$: if $(y_1,\ldots,y_{n+d+1})$
is an off diagonal prolongation of $(y_1,\ldots,y_{n+d+1})\in E_n^{(n+d+1)}$,
then for $n+1\leq j\leq n+d+1$ we have that the points $(y_1,\ldots,y_{n},
y_j)\in E_{n+1}^{(n+d+1)}$ are distinct, therefore the degree of the map is at
least $d+1$.

\subsection{Proof of theorem.} We show by induction on $i$ that for any $n$ and
 $i$ the dimension of any fiber of $F_{n+1}\rar F_n$ is at least $i+1$.
Lemma 4 shows this for $i=0$. Assume it holds true for $i-1$, let $n\geq 0$ and
let $G$ be a component of $F_{n}$, such that the fiber dimension of $F_{n+1}$
over $G$ is $i$. Applying the inductive assumption to each $F_{n+j+1}\rar
F_{n+j}$, we have that the dimension of every fiber of $F_{n+k}$ over
$F_n$ is at least $ik$. On the other hand, $F_{n+k}$ is a subscheme of the
fibered power $(F_{n+1})^k_{F_n}$, so over $G$ it has fiber dimension precisely
$ik$. Therefore
there exists a component $H_k$ of $F_{n+k}$ dominant over $G$ of fiber
dimension
$ik$, which is therefore identified as a dominant component of the fibered
power $(F_{n+1})^k_{F_n}$. By proposition 1, for $k$ sufficiently large we
have that $H_k$ is a
Lang variety. Lang's conjecture implies that $H_k(K)$ is not dense in $K$,
contradicting the definition of $F_{n+k}$. \qed

\section{A few refinements and applications in arithmetic and geometry}

\subsection{Proof of Theorem \ref{Z(X)}}
Assume that $X\rar B$ is a family of varieties of general type. By Hironaka's
desingularization theorem, we may assume that $B$ is a smooth projective
variety. Let $L$ be an
ample line bundle on $B$, let $n>>0$ be a sufficiently large integer and let
$H$ be a smooth
divisor of $L^{\otimes n}$. Let $B_1\rar B$ be the cyclic cover ramified to
order $n$ along $H$. Then by adjunction, $B_1$ is a variety of general type.
Let $X_1\rar X$ be the pullback of $X$ to $B_1$. By the main theorem (Satz III)
of \cite{viehweg}, the variety $X_1$ is of general type. Assuming the geometric
Lang conjecture, Let $Z_1(X_1)$  be the Langian exceptional set. Let
$\tilde{Z}$ be the image of $Z_1(X_1)$ in $X$. Then for any $b\in B$, we have
by definition that $Z(X_b)\subset \tilde{Z}$. \qed

It has been noted in \chm that Viehweg's work goes a long way towards proving
conjecture H. It is therefore not surprising that it may be used on occasion to
replace the assumption of conjecture H.

\subsection{Uniformity in terms of the degree of an extension}
Let $X \rar B$ be a family of GeM varieties over $K$. Assuming the conjectures,
theorem 1 gave us a uniform bound on the number of rational points over finite
extension fields in the fibers. We will now see that this in fact implies a
much stronger result, namely our theorem \ref{unideg}: the uniform bound only
depends on the degree of the field extension.

{\bf Proof of theorem \ref{unideg}:}
for $n=1$ or $2$, Let $Y_n = \Sym^d(X^n_B)$, and $Y_0 =\Sym^d(B)$.
Then we have natural maps $p_n:Y_n \rar Y_{n-1}$. Let $\Gamma$
be the branch locus of the quotient map $X^d \rar Y_1$, namely the set of
points which are fixed by some permutation. If
$P\not\in \Gamma$ then $p_2^{-1}(P)$ is a GeM variety, isomorphic over
$\overline{K}$ to the product of $d$ fibers of $X$. Denote $Y_1' = Y_1\setminus
\Gamma_1$, and $Y_2'=p_2^{-1}Y_1'$. Then $Y_2'\rar Y_1'$ is a family of GeM
varieties, and by Theorem 1 we have a bound on the cardinality of
$(Y_2')_y(K)$ uniform over $y\in Y_1'(K)$.

By induction, it suffices to bound the number of points in
$X_b(L)$ over any field $L$ of degree $d$ over $K$, which are
defined over $L$ but not over any intermediate field. If
$\sigma_1,\ldots,\sigma_d$ are the distinct embeddings of $L$ in
$\overline{K}$,
and $P\in X_b(L)$ not defined over an intermediate field, then the points
$\sigma_i(P)\in X_{\sigma_i(b)}(\sigma_i(K))\subset
X(\overline{K})$ are distinct. If $(P_1,P_2)\in X_B^2(L)$ is a
pair of such points, then the Galois orbit $\{\sigma_i(P_1,P_2),
i=1,\ldots,d\}$ is Galois stable, therefore it gives rise to a point in
$Y_2(K)$. This point has the further property that its image in $Y_1$ does not
lie in $\Gamma_1$, so it gives rise to a point in $Y'_2(K)$. The
previous paragraph shows that the number of  points on a fiber is bounded. \qed

Applying theorem 3 where $X\rar B$ is the universal family over the Hilbert
scheme of 3-canonical
curves of genus $g$ (as in \chm, \S\S 1.2), we obtain the following:

\begin{cor} Assume that the weak Lang conjecture as well as conjecture H hold.
Fix integers $d, g>1$ and a number field $K$.  Then there is a uniform bound
$N_d$ such that for
any field extension $L$ of $K$ of degree $d$ and every curve $C$ of genus $g$
over $L$ we have $\sharp C(L)<N_d$.
\end{cor}

We remark that in the cases of degrees $d\leq 3$ one does not need to assume
conjecture H: this was proven in \cite{abr}, using the fact that conjecture H
holds for families of curves or surfaces. A similar result is being
worked out by P. Pacelli  for arbitrary $d$.

Here is a special case: let $f(x)\in
{\Bbb Q}(x)$ be a polynomial of degree $>4$ with distinct complex roots. Then,
assuming the weak Lang conjecture,
the number of rational points over any quadratic field on the curve $C: y^2 =
f(x)$ is bounded uniformly. We remark that, if $\deg f>6$, this in
fact may be
deduced using a combination of \cite{chm} and a theorem of Vojta \cite{vojta}
which says that all but finitely many quadratic points on $C$ have rational $x$
coordinate. One then applies \cite{chm} which gives a uniform bound on the
rational points on the twists $ty^2 = f(x)$.

Following the suggestion of \cite{chm}, \S 6 one can apply Theorem 1 to
symmetric powers of curves. Since conjecture H is known for surfaces, one
obtains the following (stated without proof in \cite{chm}, Theorem 6.2):

\begin{cor}(Compare \chm, Theorem 6.2) Assume that the weak Lang conjecture
holds.  Fix a number field $K$.
Then there is a uniform bound $N$ for the number of quadratic
points on any nonhyperelliptic, non-bielliptic curve $C$ of genus $g$
over $K$.
\end{cor}

Similarly, it was shown in \cite{ah}, lemma 1 that if the gonality of a curve
$C$ is $>2d$ then $\Sym^d(C)$ is GeM. Recall that a closed point $P$ on $C$ is
said to be
of degree $d$ over $K$ if $[K(P):K]=d$. We deduce the following:

\begin{cor} Assume that the weak Lang conjecture holds.  Fix a number field $K$
and an integer $d$.
Then there is a uniform bound $N$ for the number of
points of degree $d$ over $K$ on any curve $C$ of genus $g$ and gonality $>2d$
over $K$.
\end{cor}

\subsection{The geometric case}
One can use the same methods using Lang's conjecture for function fields of
characteristic 0, say over $\Bbb C$.  Given a  fibration
$X\rar B$ where the generic fiber is a variety of general type, a rational
point
$s\in X(K_B)$ over the function field  of $B$
is called {\em constant} if $X$ is birational to a product $X_0\times B$ and
$s$ corresponds to a point on $X_0$.  Lang's conjecture for function fields
says
that the non-constant points are not Zariski dense.

In this section we will restrict attention to the case where the base is the
projective line $\BP$. We will only assume the following statement: if $X$ is
a variety of general type, then the rational curves in $X$ are not Zariski
dense. It is easy to see that this statement in fact follows from  the
geometric Lang conjecture, as well as from Lang's conjecture for function
fields.

We would like to apply this conjecture to obtain geometric uniformity results.
One has to be careful here, since the conjecture does not
apply to Lang varieties, and one has to use a variety of general type directly.

As stated in the introduction, if $X\rar B$ is a family of curves of genus $>1$
the appropriate variety $W$  of general type dominated by $X^n_B$ is identified
in \chm as the  image
$B_n\subset{\MM}_{g,n}$ of $X^n_B$ by the moduli map. We use this in the proof
of the following proposition:

\begin{prp}\label{unigeom} Assume that Lang's conjecture for function fields
holds.  Fix
an integer $g>1$.
Then there is a  bound $N$ such that for any generically smooth family of
curves $C\rar \Bbb{P}^1$ of genus $g$ there are at most $N$ non-constant
sections $s:\Bbb{P}^1\rar C$.
\end{prp}

{\bf Proof.} First note that if $s:\Bbb{P}^1\rar C$ is a nonconstant section
whose image in $ \MM_{g,1}$ is  a point, then $s$ becomes a
constant section after
a finite base change $D\rar \Bbb{P}^1$. This implies that $s$ is  fixed by
 a nontrivial automorphism of $C$, and the number of such points is bounded
in terms of $g$.
Therefore it suffices to bound the number of sections whose image in
$\MM_{g,1}$
is non-constant. We will call such sections {\em strictly non-constant}.

 Let $B_0\subset \MM_g$ be a closed subvariety, and choose $n$  such
that $B_n\subset \MM_{g,n}$ is of general type. If a family $C\rar \BP$ has
moduli in $B_0$, then for any $n$-tuple of strictly
non-constant
sections $s_i:\BP\rar C$, we obtain a non-constant rational map  $\BP \rar
B_n$. Let ${F} \subset
B_n$ be the Zariski closure of the images of the collection of non-constant
rational maps obtained this way.

Since $B_n$ is of general type, Lang's
conjecture implies that $F\neq B_n$. Applying
lemma 1.1 of \chm we obtain that there is an closed subset set $F_0\subset B_0$
and an
integer $N$ such that, given a family of curve $C\rar\BP$ such that the
rational image
of $\BP$ in $\MM_g$ lies in $B_0$ but not in  $F_0$, there are at most $N$
non-constant sections of $C$. Noetherian induction gives the theorem. \qed

Choosing a coordinate $t$ on  $ {\Bbb P}^1$  we can pull back the curve $C$
along the map $\BP\rar \BP$ obtained by
taking $n$-th roots of $t$. Let ${\Bbb C}(t^{1/\infty}) = {\Bbb C}(\{t^{1/n},
n\geq 1\})$, the field obtained by adjoining all roots of $t$. If one restricts
attention to non-isotrivial curves, one obtains the following amusing result
(suggested to the author by Felipe Voloch):

\begin{cor} Assume that the Lang conjecture for function fields holds.  Fix an
integer $g>1$.
Then there is a  bound $N$ such that for any smooth
nonisotrivial
curve $C$ over ${\Bbb C}(t)$ of genus $g$ there are at most $N$ points in
$C({\Bbb C}(t^{1/\infty}))$.
\end{cor}

One can also try to prove uniformity results analogous to theorem 3. Using the
results in \cite{abr} we can refine proposition \ref{unigeom} and obtain
theorem \ref{unigon}.

{\bf Proof of theorem \ref{unigon}.} The proof is a slight modification of the
theorem of \cite{abr}, keeping track of the dominant map to a variety of
general type.

 As in the proof of theorem 3, it suffices
to look at sections $s:D\rar C$ which are not pullbacks of  sections of a
family
over $\BP$.

 In an analogous way to the
proof of theorem \ref{unif}, we say that an
$n$-tuple of distinct, strictly  non-constant sections is $m$-prolongable if
it may be prolonged to an $m$-tuple of distinct, strictly non-constant
sections, none of which being the pullback from a family over $\BP$.  Any
$n$-tuple of distinct
sections  $s_i:D\rar C$ over a hyperelliptic curve $D$ gives rise to a rational
map
$\BP\rar  \Sym^2(\MM_{g,n})$. We define $F_n^{(m)}$ to be the
closure in $\Sym^2(\MM_{g,n})$ of the images of $m$-prolongable sections, and
$F_n = \cap_{m>n} F_n^{(m)}$.

As in Lemma 1, we have that the relative
dimension of
any fiber of $F_{n+1}\rar F_n$ is positive. We have two cases to consider:
either for high $n$ there is a component of $F_{n+1}$ having fiber dimension 1
over $F_n$, or for all $n$ the fiber dimension is everywhere 2.

In case the fiber dimension is 1,
we will see that there is a
component  of $F_{n+k}$ which is a variety of general type.
Assuming Lang's conjecture for function fields this contradicts the fact that
the images of non-constant sections are dense. Fix a
 general fiber $f$ of $F_{n+1}$ over $F_{n}$. The curve $f$ lies inside a
surface isomorphic to the product of two curves
$C_{b_1} \times C_{b_2}$. By the definition of $m$-prolongable sections, and
analogously to lemma 1, we obtain that there is a component $f'$ of $f$ which
maps surjectively to both $C_{b_1}$ and $C_{b_2}$. Therefore as either $b_1$ or
$b_2$ moves in $B_0$, the moduli of $f'$ move as well.

Let $F'$ be a component
of
$F_{n+1}$ whose fibers have this property. If we follow the proof of
proposition 1 and use the moduli description of the dominant map to a variety
of general type $m:(F')^k_{F_n}\rar W$, we see that if $E$ is a general curve
in
$(F')^k_{F_n}$ lying in a fiber of $m$, then $E$ projects to a point
in $B_0$; moreover, by the definition of prolongable points, $E$ projects to an
off diagonal point in some $(F')^l_{F_n}$. But the fibers over off-diagonal
points are GeM varieties, therefore the general fiber of the map $m$ is of
general type. By the main theorem of \cite{viehweg},
$(F')^k_{F_n}$ is itself a variety of general type, and therefore $F_{n+k}$ has
a component of general type, contradicting Lang's conjecture.

In case of fiber dimension 2, we
use
proposition 1 of \cite{abr}: let $B\subset \Sym^2(\MM_g)$. Then for high $n$,
the inverse image
$B_n\subset \Sym^2(\MM_{g,n})$ of $B$ is a variety of general type. Since the
images of non-constant sections are dense in $F_n$, this again
contradicts Lang's conjecture. \qed

 If one restricts attentions to trivial fibrations, one obtains as an immediate
corollary:

\begin{cor} Assume that the Lang conjecture for function fields holds.  Fix an
integer $g>1$.
Then there is an integer $N$ such that for any
curve $C$ of  genus $g$ and any hyperelliptic curve $D$
there are at most $N$ non-trivial morphisms $f:D\rar C$.
\end{cor}

It should be noted that the theory of Hilbert schemes gives the existence of a
bound depending on the genus of $D$, which is clearly not as strong. As in the
arithmetic case, I expect that work in progress of Pacelli should generalize
these results to the case where $D$ is $d$-gonal, for fixed $d$.

\today\\
Dan Abramovich\\
Department of Mathematics, Boston University\\
111 Cummington Street, Boston, MA 02215 \\
abrmovic@@math.bu.edu \\
http://math.bu.edu/INDIVIDUAL/abrmovic


\begin{thebibliography}{HHHHHHH}

 \bibitem[$\aleph$]{abr} D. Abramovich: {\em Uniformit\'e des points rationnels
des courbes alg\'ebriques sur les extensions quadratiques et cubiques},
 C.R. Acad. Sc. France, to appear.
{\tt http://math.bu.edu/INDIVIDUAL/abrmovic/quadf.ps }

\bibitem[AH]{ah} D. Abramovich, J. Harris: {\em
 Abelian varieties and curves in  $W_d(C)$.}
Compositio Math. 78 (1991) p. 227-238.

 \bibitem[CHarM]{chm}  L. Caporaso, J. Harris, B. Mazur: {\em Uniformity of
   rational points,}  J. Amer. Math. Soc., to appear. Temporarily available on
{\tt http://math.bu.edu/INDIVIDUAL/abrmovic/chmlong1.ps }

 \bibitem[Has]{hassett} B. Hassett: {\em Correlation for surfaces of general
type},
preprint, 1995.

 \bibitem[LangAMS]{langbul}  S. Lang, {\it Hyperbolic diophantine analysis.}
   Bull. A.M.S. 14 (1986) p. 159-205.

 \bibitem[LangIII]{lang3}  S. Lang, {\it Number Theory III: Diophantine
Geometry.}
Encyclopedia of Mathematical Sciences vol. 60. Springer Berlin - New York 1991.


\bibitem[Vie]{viehweg} E. Viehweg:  {\em Die Additivit\"at der Kodaira
Dimension f\"ur projektive Fasserr\"aume \"uber Variet\"aten des allgemeinen
Typs,} Jour. reine und angew. Math. 330 (1982), 132-142.


\bibitem[Voj]{vojta} P. Vojta: {Arithmetic discriminants and quadratic points
in
curves.} In: van der Geer, Oort, Steenbrink, eds., {\em Arithmetic algebraic
geometry Texel 19889}, Progress in Mathematics 89, Birkh\"auser, Boston 1991. P
359-376.

\bibitem[Vol]{voloch} J. F. Voloch:
{\em Diophantine geometry in characteristic $p$: a survey}, preprint.
{\tt ftp://math.utexas.edu/pub/papers/voloch }


\end{thebibliography}
\end{document}